\documentclass[12pt]{article}
\usepackage{amsmath,euscript}
\usepackage{amssymb}
\usepackage{amsthm}
\usepackage{enumerate}
\usepackage{amsfonts}
\usepackage{comment}
\usepackage{colortbl}
\usepackage{a4wide}
\usepackage{amssymb,amsmath,array}
\usepackage{amssymb}
\usepackage{amsthm}
\usepackage{enumerate}
\usepackage{amsfonts}
\usepackage{comment}
\usepackage{colortbl}
\usepackage{a4wide}
\usepackage{amssymb,amsmath,array}

\newtheorem{thm}{Theorem}[section]
\newtheorem{lem}[thm]{Lemma}

\newtheorem{conj}[thm]{Conjecture}
\theoremstyle{defn}

\newtheorem{rem}{Remark}

\begin{document}
\bibliographystyle{plain}
\begin{center}
{\bf \large Division algebras of Gelfand-Kirillov transcendence degree $2$
}\\
\vspace{7 mm} \centerline{Jason P. Bell\footnote{The author thanks NSERC for its generous support.}}

\centerline{Department of Mathematics}
\centerline{
Simon Fraser
University}
\centerline{
8888 University Dr.}
\centerline{Burnaby, BC V5A 1S6.}
\centerline{ CANADA}
\centerline{\tt jpb@math.sfu.ca}
\centerline{AMS Subject Classification: 16P90}
\centerline{Keywords: division algebras, subfields, GK dimension, growth}
\end{center}

\begin{abstract} Let $A$ be a finitely generated $K$-algebra that is a domain of GK dimension less than $3$, and let $Q(A)$ denote the quotient division algebra of $A$.  We show that if $D$ is a division subalgebra of $Q(A)$ of GK dimension at least $2$ then $Q(A)$ is finite dimensional as a left $D$-vector space.  We use this to show that if $A$ is a finitely generated domain of GK dimension less than $3$ over an algebraically closed field $K$ then any division subalgebra $D$ of $Q(A)$ is either a finitely generated field extension of $K$ of transcendence degree at most one, or $Q(A)$ is finite dimensional as a left $D$-vector space.
\end{abstract}

\section{Introduction}
Given a finitely generated algebra $A$ over a field $K$, the \emph{Gelfand-Kirillov dimension} (GK dimension, for short) of $A$ is defined to be:
\[ {\rm GKdim}(A) \ = \ \limsup_{n\rightarrow\infty} \frac{\log\, {\rm dim} (V^n)}{\log \, n}, \]
where $V$ is a finite dimensional $K$-vector subspace of $A$ which contains $1$ and generates $A$ as a $K$-algebra.  We note that this definition is independent of the choice of vector space $V$ with the above properties.  Gelfand-Kirillov dimension should be viewed as a noncommutative analogue of Krull dimension; indeed, ${\rm GKdim}(A)={\rm Kdim}(A)$ when $A$ is a finitely generated commutative $K$-algebra.  As a result, GK dimension has been used to obtain noncommutative analogues of results from classical algebraic geometry \cite{ArtSt, StV}.
We refer the reader to Krause and Lenagan \cite{KL} for the basic facts about GK dimension.  

Just as commutative domains have a field of fractions, finitely generated noncommutative domains of finite GK dimension are Goldie rings \cite[p. 46]{KL}, and hence have a quotient division algebra.  Given a finitely generated domain $A$ of finite GK dimension, we let $Q(A)$ denote the quotient division algebra of $Q(A)$.  One of the goals in noncommutative algebra is to understand domains of GK dimension $2$; Artin and Stafford \cite{ArtSt} gave a concrete description of graded domains of GK dimension $2$, but the ungraded case is not well-understood.  Just as one has a birational classification of projective surfaces in classical algebraic geometry, one would like to obtain a ``birational'' description of algebras of GK dimension $2$; that is, to have a concrete description of the quotient division algebras of the finitely generated domains of GK dimension $2$ over an algebraically closed field.  Artin \cite{Ar} has a proposed classification, however, not much progress has been made towards obtaining this classification.  

Working towards a birational classification of domains of GK dimension $2$, one seeks invariants which can be used.  Several authors have proposed a noncommutative analogue of transcendence degree that could be useful in attaining this goal \cite{GK, Z, YZ}.  In general, one expects quotient division algebras of domains of GK dimension $d$ to share some properties with fields of transcendence degree $d$.  In particular, Zhang \cite[Conjecture 8.4]{Z} has made the following conjecture.
\begin{conj} Let $K$ be a field and let $A$ be a finitely generated $K$-algebra that is a domain of GK dimension $d$.  Suppose
$$K=D_0 \subseteq D_1 \subseteq \cdots \subseteq D_m \subseteq Q(A)$$ is a chain of division subalgebras of $Q(A)$ such that:
\begin{enumerate}
\item $D_i$ is infinite dimensional as a left $D_{i-1}$-vector space for $i=1,2,\ldots ,m$;
\item $D_i$ is finitely generated as a division algebra.
\end{enumerate}
Then $m\le d$.
\label{conj: 1}
\end{conj}
We note that this conjecture is highly non-trivial.  Indeed, if one replaces the $D_i$ with subalgebras and  then insists that $D_{i+1}$ is a free left $D_i$-module of infinite rank, then the conjecture is not true in general.  For example, the Weyl algebra, $K\{x,y\}/(xy-yx-1)$ is a domain of GK dimension $2$, and Makar-Limanov \cite{ML} has shown that its quotient division algebra contains a free algebra on infinitely many variables if $K$ has characteristic $0$.  The fact that GK dimension is poorly behaved under localization creates the need for new invariants to study quotient division algebras.

Our main result is to prove Conjecture \ref{conj: 1} for algebras of GK dimension strictly less than three.

\begin{thm} Let $K$ be a field and let $A$ be a finitely generated $K$-algebra of GK dimension strictly less than $3$ that is a domain and let $Q(A)$ denote the quotient division algebra of $A$.   Suppose $$K \ =\ D_0\ \subseteq\  D_1\  \subseteq\  D_2\  \subseteq\  \cdots \subseteq D_m\ \subseteq \  Q(A)$$ is a chain of division subalgebras of $Q(A)$ such that:
\begin{enumerate}
\item $D_i$ is infinite dimensional as a left $D_{i-1}$-vector space for $i=1,2,\ldots ,m$;
\item $D_i$ is finitely generated as a division algebra for $i=1,2,\ldots m$.
\end{enumerate}
Then $m\le 2$.
\label{thm: 1}
\end{thm}
In fact one can make the following, even stronger, statement.
\begin{thm}  Let $K$ be a field and let $A$ be a finitely generated $K$-algebra of GK dimension strictly less than $3$ that is a domain and let $Q(A)$ denote the quotient division algebra of $A$.  If $D$ is a division subalgebra of $Q(A)$ and $D$ has GK dimension at least $2$, then $Q(A)$ is finite dimensional as a left $D$-vector space.\label{thm: main1}
\end{thm}
We note that this theorem is closely related to Smoktunowicz's gap theorem, which states that a connected finitely graded domain cannot have GK dimension strictly between $2$ and $3$.  We note that if $A$ is a connected finitely graded domain of GK dimension less than $3$, then $A$ has a graded quotient $D_0[x,x^{-1};\sigma]$ and by Theorem \ref{thm: main1}, $D_0$ must have GK dimension at most $1$.  Thus $A$ has GK dimension exactly $2$ by a result of Artin and Stafford  \cite[Theorem 0.5, p. 242]{ArtSt}.

An interesting consequence of Theorem \ref{thm: 1} is the following.

\begin{thm} Let $K$ be an algebraically closed field and let $A$ be a finitely generated $K$-algebra that is a domain of GK dimension $< 3$.  If $D$ is a division subalgebra of $Q(A)$ then either $D$ is commutative or $Q(A)$ is finite dimensional as a left $D$-vector space; moreover, if $D$ is commutative and has transcendence degree at least $2$ over $K$ then $Q(A)$ is finite dimensional over its centre.
\label{thm: main2}
\end{thm}

We use some of the ideas of Smoktunowicz \cite{Sm} in studying division subalgebras of quotient division algebras of domains of GK dimension less than $3$.  
 Smoktunowicz's key idea in studying chains of division algebras
$$ Z(D) =F \subseteq E\subseteq D,$$ lies in the observation that if  $S$ is a subset of $E$ that is linearly independent over $F$ and $T$ is a subset of $D$ that is right-linearly independent over $E$ then the set
$\{ ts~|~s\in S, t\in T\}$ is linearly independent over $F$.  Using this observation, in addition to some estimates, we are able to prove Theorems \ref{thm: 1}, \ref{thm: main1}, and \ref{thm: main2}.

Theorems \ref{thm: 1} and \ref{thm: main2} are proved in Section 2.  In Section 3, we give some concluding remarks and make some conjectures.

\section{Quotient division algebras of GK 2 domains}
In this section we prove theorem \ref{thm: 1}.  To do this, we need a few estimates.
We begin with an argument that is based on an argument of Smoktunowicz \cite{Sm}.
\begin{lem} Let $K$ be a field and let $D$ be a division ring that is a $k$-algebra.  Suppose that $V$ is a finite dimensional vector subspace of $D$ such that
$1\in V$ and ${\rm dim}(V^n) > d\, {\rm dim}(V^{n+1}/V^{n})$.  If $s\in V$ and $x_1,\ldots ,x_d\in s^{-1}V$, then there is nonzero
$u\in V^n$ such that $ux_1\cdots x_i\in V^n$ for $1\le i\le d$.   
\label{lem: 1}
\end{lem}
\noindent {\bf Proof.} Write $x_i=s^{-1}a_i$ with $a_i\in V$ for $1\le i\le d$.
 Let $W_{0,n}=V^n$ and
for $1\le i\le d$, 
let $W_{i,n} = \{ r\in W_{i-1,n}~|~rx_1\cdots x_i \in V^n\}$, and let
$$U_{i,n} = W_{i,n}x_1\cdots x_i.$$
Observe that
$${\rm dim} \,W_{i,n} \ = \ {\rm dim} \, W_{i-1,n}\cap V^n(x_1\cdots x_i)^{-1} \ = \  {\rm dim} \, U_{i-1,n}\cap V^n x_i^{-1} \ =  \ {\rm dim} \,U_{i-1,n}\cap V^{n} sa_i^{-1}.$$
Thus
\begin{eqnarray*}
{\rm dim}(W_{i,n})  &\ge & {\rm dim}\, U_{i-1,n}a_i \cap V^{n} s \\
&=& {\rm dim}(U_{i-1,n}a_i) + {\rm dim}(V^{n}s) - {\rm dim}(U_{i-1,n}a_i\cup V^{n}s) \\
&=& {\rm dim}(W_{i-1,n}) + {\rm dim}(V^{n}) - {\rm dim}(V^{n+1}) \\
&=& {\rm dim}(W_{i-1,n}) - {\rm dim}(V^{n+1}/V^{n}).
\end{eqnarray*}
Hence
$${\rm dim}(W_{i,n}) - {\rm dim}(W_{i-1,n}) \ge - {\rm dim}(V^{n+1}/V^{n}).$$
It follows that
\begin{eqnarray*}
{\rm dim}(W_{d,n}) &=& 
 {\rm dim}(W_{0,n}) - \sum_{i=0}^{d-1} {\rm dim}(W_{i,n}) - {\rm dim}(W_{i+1,n}) \\
&\ge & {\rm dim}(V^n) - \sum_{i=0}^{d-1} {\rm dim}(V^n/V^{n-1}) \\
&=& {\rm dim}(V^n) - d\,  {\rm dim}(V^{n+1}/V^{n}).\end{eqnarray*}
The result follows. \qed
\begin{rem} Suppose $A$ is a finitely generated Goldie domain over a field $K$ and $E$ is a division subalgebra of $Q(A)$.  If $Q(A)$ is infinite dimensional over $E$ as a right $E$-vector space then
$AE$ is infinite dimensional as a right $E$-vector space.\label{rem: 1}
\end{rem}
\noindent {\bf Proof.} Suppose that ${\rm dim}_E(AE)=m$.  By assumption there exist $\alpha_1,\ldots ,\alpha_{m+1}\in Q(A)$ that are right-linearly independent over $E$.  There exists some $b\in A$ such that $a_i:=b\alpha_i\in A$ for $1\le i\le m+1$.  Then by construction, $a_1E+\cdots +a_{m+1}E$ is direct, contradicting the fact that ${\rm dim}_E(AE)$ is finite dimensional. \qed

\vskip 2mm
\noindent The following lemma is used to construct a large linearly independent set inside a division algebra.  This will eventually be combined with Lemma \ref{lem: 1} to get lower bounds on the GK dimension of a domain.  
\begin{lem} Suppose that $D$ is a finitely generated division algebra with centre $K$ and $E$ is a division subalgebra generated (as a division algebra) by a finite dimensional $K$-vector space $S$.  If $E$ has GK dimension at least $2$ and $K$ has infinite transcendence degree over its prime subfield, then there exist elements $a_i,b_i\in S$ for $i\ge 1$ such that the set of elements 
$\{a_1\cdots a_i b_j \cdots b_1~|~i,j\ge 1\}$ are linearly independent over $K$.
\label{lem: newsmok}
\end{lem}
\noindent {\bf Proof.} Cf. Smoktunowicz \cite[Lemma 2]{Sm3}.\qed

\begin{lem} Let $A$ be a finitely generated $K$-algebra of finite GK dimension and let $V$ be a finite dimensional generating subspace for $A$ that contains $1$.  Then there is a positive constant $C$ such that
$${\rm dim}(V^n) \ > \ C n\, {\rm dim}(V^{n+1}/V^{n})$$ for infinitely many $n$. \label{lem: 2}
\end{lem}
\noindent {\bf Proof.} Since $A$ has finite GK dimension, there is some $d>1$ and some constant $C_0>0$ such that
$${\rm dim}(V^n) \ < \ C_0 n^d \qquad {\rm for~all}~n>0.$$
Suppose that the result does not hold.  Then there exists some $N$ such that
$${\rm dim}(V^n) \ < \ \frac{1}{4^d} n\, {\rm dim}(V^{n+1}/V^{n})$$ for all $n\ge N$.
Thus for $n\ge N$ we have
\begin{eqnarray*}
{\rm dim}(V^{2n}) &=& {\rm dim}(V^N) + \sum_{i=N}^{2n-1} {\rm dim}(V^{i+1}/V^i) \\
&\ge & {\rm dim}(V^N) + 4^d \sum_{i=N}^{2n-1} \frac{1}{i} {\rm dim}(V^{i})\\
&\ge & 4^d \sum_{i=n}^{2n-1} \frac{1}{i} {\rm dim}(V^{i+1}) \\
&\ge & 4^d {\rm dim}(V^n)(1/n+\cdots + 1/(2n-1)) \\
&\ge & 2^{2d-1} {\rm dim}(V^n). \end{eqnarray*}
Thus by induction we see
$$C_0(2^{\ell}N)^d \ge {\rm dim}(V^{2^{\ell}N}) \ \ge \ {\rm dim}(V^N)2^{(2d-1)\ell}$$
for all $\ell>0$.  
Taking the $\ell$'th root and taking the limsup as $\ell$ tends to infinity, we see
$$2^d \ge \ 2^{2d-1},$$ a contradiction since we picked $d>1$. 
\qed
\vskip 2mm
We are now ready to prove our main results.
\vskip 2mm
\noindent {\bf Proof of Theorem \ref{thm: main1}.} We first note that $A\otimes_K K(t_1,t_2,\ldots)$ is a domain since it is a localization of $A[t_1,t_2,\ldots]$.  Thus it is no loss of generality to assume that $K$ has infinite transcendence degree over its prime subfield.  Let $S$ be a finite dimensional $K$-vector space that generates $D$ as a division algebra.  Since $D$ has GK dimension at least $2$, 
by Lemma \ref{lem: newsmok}
there exist $a_i,b_i\in S$ such that
$$\{ a_1\cdots a_ib_j\cdots b_1~|~i,j\ge 1\}$$ is a linearly independent subset of $D$ over 
$K$.  

Let $V$ be a finite dimensional $K$-vector subspace of $A$ such that:
\begin{enumerate}
\item there are $b,b '\in V$ such that $bS, Sb' \subseteq V$;
\item $1\in V$;
\item $V$ generates $A$ as a $K$-algebra.
\end{enumerate}
 Suppose $A$ has GK dimension less than $3$.  By Lemma \ref{lem: 2}, there exists a constant $C>0$ such that
$${\rm dim}(V^n)  >  C n\, {\rm dim}(V^{n+1}/V^{n})$$ for infinitely many $n$.  Let 
\begin{equation}
T \ = \ \left\{ n~:~{\rm dim}(V^n) \ > \ C n\, {\rm dim}(V^{n+1}/V^{n})\right\}.
\end{equation}
Let $n\in T$.
By Lemma \ref{lem: 1}, there exists some $r\in V^{n}$ such that
$ra_1\cdots a_i\in V^n$ for $1\le i<Cn$.
Similarly, there exists some $s\in V^n$ such that
$b_1\cdots b_js\in V^n$ for $1\le j<Cn$.
By assumption, 
$$\{ a_1\cdots a_ib_j\cdots b_1~|~1\le i,j\le n\}$$ is a subset of $D$ that is linearly independent over $K$. 
Since $Q(A)$ is left infinite dimensional over $D$, we see that there exist $u_1,\ldots ,u_n\in V^n$ such that
$(s^{-1}Ds)u_1+\cdots + (s^{-1}Ds)u_n$ is direct by Remark \ref{rem: 1}.
Thus
\[ \left\{ s^{-1} a_1\cdots a_j b_1\cdots b_k s u_i ~|~1\le i\le n, 1\le j,k< Cn\right\} \] is a $K$-linearly independent subset of $Q(A)$.
Hence
\[ \left\{  r a_1\cdots a_j b_1\cdots b_k s u_i ~|~1\le i\le n, 1\le j,k < Cn\right\} \] is $K$-linearly independent.
Since $u_i\in V^n$, $ra_1\cdots a_j\in V^n$, and $b_1\cdots b_k s\in V^n$ for $1\le i,j,k\le n$, we see that
\[ \left\{ r a_1\cdots a_j b_1\cdots b_k su_i ~|~1\le i\le n, 1\le j,k<Cn \right\} \ \subseteq V^{3n}.\]
It follows that
\[ {\rm dim}(V^{3n}) \ \ge \ (Cn-1)^2 n \qquad {\rm for}~n\in T\] and thus
$$\limsup\, \log ({\rm dim}(V^n))/\log\, n \ \ge \ 3.$$  The result now follows. \qed   
\vskip 2mm

\noindent {\bf Proof of Theorem \ref{thm: 1}.}  Suppose we have a chain of division subalgebras
$$D_0\subseteq D_1\subseteq D_2\subseteq D_2\subseteq Q(A)$$ satisfying the hypotheses in the statement of the theorem.  We show that $A$ must have GK dimension at least $3$.  Since $D_1$ is left infinite dimensional and finitely generated over $K$ as a division algebra, it must have GK dimension at least $1$.  Similarly, since $D_2$ is finitely generated as a division algebra and is left infinite dimensional over $D_1$ it must have GK dimension at least $2$.   To see this, let $S_i$ be a finite set containing $1$ generating $D_i$ as a division algebra.
Since $D_{2}$ is left infinite dimensional over $D_{1}$, there exist
$a_1,\ldots ,a_n$ in $S_2^n$ that are left linearly independent over $D_1$ by Remark \ref{rem: 1}.  Similarly, since $D_1$ has GK dimension at least $1$, there exist $b_1,\ldots ,b_n\in S_1^n$ that are right linearly independent over $K$.  Consequently, the set 
$$\{b_ja_i~|~1\le i,j\le n\}\subseteq (S_1\cup S_2)^n$$ is linearly independent over $K$.  It follows that the GK dimension of $D_2$ is at least $2$.  Since $Q(A)$ is infinite dimensional over $D_2$ as a left $D_2$-vector space, we see that $A$ has GK dimension at least $3$ by Theorem \ref{thm: main1}. \qed   

\vskip 2mm
We prove Theorem \ref{thm: main2}.
\vskip 2mm
\noindent {\bf Proof of Theorem \ref{thm: main2}.} Let $K$ be a subfield of $Q(A)$.  If $K$ has transcendence degree at least $2$, then $Q(A)$ must be finite dimensional as a left $K$-vector space by Theorem \ref{thm: main1}.  But this gives that $A$ is PI, and hence $Q(A)$ is finite dimensional over its centre.  Let $D$ be a division subalgera of $A$.  If $Q(A)$ is not finite dimensional as a left $D$-vector space, then $D$ must be of GK dimension $1$ by Theorem \ref{thm: main1} and hence it is locally PI by the Small-Warfield theorem \cite{SW}.  
We claim that $D$ is commutative.   Suppose that $D$ is noncommutative and pick $x$ and $y$ in $D$ that do not commute with each other.  Let $D_0$ denote the division algebra generated by $x$ and $y$.  Since $D_0$ has GK dimension $1$, $D_0$ is finite dimensional over $K(x)$ \cite[Theorem 4.12]{KL} and hence it is PI.  This means that $D_0$ is finite dimensional over its centre, and since $D_0$ has GK dimension $1$, its centre is a finitely generated field extension of $K$ of transcendence degree $1$.  Hence $D_0$ is is commutative by Tsen's theorem.   \qed


 
\section{Concluding remarks}
The techniques used in this paper are not able to obtain a proof of Conjecture \ref{conj: 1} because two-sided estimates are used.  It is probably necessary to develop new, stronger invariants that require only one-sided estimates.  Complicating matters is the fact that there exist a division algebra $D$ with a division subalgebra $E$ such that $D$ is finite dimensional as a left $E$-vector space and is infinite dimensional as a right $E$-vector space; these examples do not have finite stratiform length, however\cite{Sc}. Working towards Conjecture \ref{conj: 1}, we make the following, theoretically easier, conjecture.

\begin{conj} Let $K$ be a field and let $A$ be a finitely generated $K$-algebra that is a domain of finite GK dimension.  Suppose that $D$ is a division subalgebra of $Q(A)$ such that $D\cong Q(A)$.  Then $Q(A)$ is finite dimensional as a left $D$-vector space.
\end{conj}
We note that if $K$ is a field of characteristic zero and $A$ is the $K$-algebra formed by taking the enveloping algebra of the positive part of the Witt Lie algebra, then $A$ is a finitely generated Goldie domain of infinite GK dimension and $Q(A)$ has a division subalgebra $D$ such that $D\cong Q(A)$ and $Q(A)$ is infinite dimensional as a left $D$-vector space.

Theorem \ref{thm: main2} is related to a result of Smoktunowicz \cite{Sm2} and to a result of Small and the author \cite{BS}.  Smoktunowicz \cite{Sm2} shows that the centralizer of a non-algebraic element of a domain of  quadratic growth over a finite field is a PI domain.  One can eliminate the hypothesis that $K$ be algebraically closed in the statement of Theorem \ref{thm: main2}, but the conclusion of the statement will then be that a division subalgebra $D$ of $Q(A)$ is either locally PI (that is, a union of finitely generated division algebras, each finite dimension over their centres) or $Q(A)$ is finite dimensional as a left $D$-vector space.  In the case of an algebraically closed field, Tsen's theorem allows us to replace PI with the stronger condition of being commutative.  The author and Small \cite{BS} consider centralizers in non-PI domains of GK dimension $2$ and show that over an algebraically closed field the centralizer of a non-scalar element is commutative.  Smoktunowicz \cite{Sm4} proves an analogue of this result for quotient division algebras, without the hypothesis that the field be algebraically closed.  

We make the remark that the results we give can all be expressed more generally using prime Goldie rings instead of domains.  

\section{Acknowledgments}
I thank Agata Smoktunowicz for many helpful comments.

\end{document}